\documentclass[titlepage,a4paper]{article}
\usepackage{times}
\usepackage{mathptm}
\usepackage{amsfonts}
\usepackage{amsmath}
\usepackage{amssymb}
\usepackage{amsthm}
\usepackage[dvips]{graphicx} 
\newtheorem{thm}{Theorem}[section]
\newtheorem{lemma}{Lemma}[section]
\newtheorem{prop}{Proposition}[section]
\newtheorem{cor}{Corollary}[section]

\def \mr {{\mathbb R}}

\def \lap {\Delta}

\def \p {\partial}

\def \rao#1 {\frac{\p}{\p #1} #1}

\def \la {\langle}
\def \ra {\rangle}
\newcommand{\contr}[2]{
          {#1}\,
           \reflectbox{\rotatebox[origin=c]{180}{\mbox{$\neg$}}}
          \,{#2}}%

\newcommand{\psDO}[0]{{\Psi\!\mathrm{DO}}}%
\newcommand{\psDOcl}[1]{{\Psi\!\mathrm{DO}^{#1}_{\mathrm{cl}}}}%
\numberwithin{equation}{section}

\begin{document}
\title{An inverse boundary value problem for harmonic differential forms}
   \author{M. S. Joshi
\\Royal Bank of Scotland Group
  Risk  ,
\\135 Bishop's Gate,
\\ London EC2M 3RU,
\\E-mail: {\tt \small mark\_joshi@yahoo.com}
\and
W.R.B. Lionheart
\\Department of Mathematics,
\\UMIST,
\\P.O. Box 88,
\\Manchester M60 1QD,
\\E-mail: {\tt \small Bill.Lionheart@umist.ac.uk} }
\date{24th June 2002}
\maketitle

\begin{abstract}
We show that the full symbol of the Dirichlet to Neumann map of the
$k$-form Laplace's equation on a Riemannian manifold (of dimension
greater than 2) with boundary determines the full Taylor series, at
the boundary, of the
metric. This extends the result of Lee and Uhlmann for
the case $k=0$. The proof avoids the computation of the full symbol by using the
calculus of pseudo-differential operators parametrized by a boundary
normal coordinate and  recursively calculating the principal symbol of
the difference of boundary operators. 
\end{abstract}
\section{Introduction}
While there many results on the
uniqueness of recovery of the coefficients of an elliptic partial
differential equation  from boundary data in the case of a single
partial differential equation, there are few results for systems of
PDEs. One might expect that the complete boundary data for a system
might be sufficient to recover multiple coefficients and yet the
results to date have been in essentially scalar cases. Lee and Uhlmann
showed that the full Taylor series, at the boundary, of a metric  can be
obtained from the total symbol of the Dirichlet to Neumann map of the scalar
Laplace's equation. One might expect to obtain at least the same information
from the Dirichlet to Neumann map associated with the Laplacian
operator on $k$-forms. The 1-form, or vector Laplacian on 3-manifolds
being the example with the most obvious applications. We will show
here that the full symbol of the $k$-form Laplacian does indeed
determine the Taylor series of the metric at the boundary, and hence
under suitable assumptions an analytic metric can be recovered from
this data. The method used, in common with \cite{jmcd}, avoids the
computation of the full symbol of the Dirichlet to Neumann map.
Rather, by using the calculus of pseudo-differential operators
parametrized by a boundary normal coordinate, we recursively calculate
the principal symbol of the difference of boundary operators, checking
that it vanishes to a suitable order. As well as being an natural
extension of the scalar case  and interesting in its own right, we
hope that this paper will stimulate the use of the technique in other
inverse boundary value problems for elliptic systems of equations of
interest to applications.

The context for all will be a smooth compact orientable
manifold with boundary $M$, equipped with a Riemannian metric $g$. We
also assume $\dim M= n>2$. The metric tensor induces a volume form
$\mu \in \Omega^n(M)$ and Hodge star isomorphism
$*:\Omega^k(M)\rightarrow \Omega^{n-k}(M)$ is defined by the property
\begin{equation}\label{equation:startensor} *\omega \wedge \omega =
g( \omega, \omega) \mu,
\end{equation}
where the action of the metric is extended naturally to act on
$k$-forms. We can consider the Hodge star on $k$-forms as a
contraction of the tensor $g^{\sharp\,\otimes k}\otimes \mu$. Here
$g^\sharp$ is the covariant metric tensor.

The total symbol of an operator $P$ on functions on $\mathbb{R}^n$ is 
$$p(x,\xi)=e^{-ix\cdot\xi}P(e^{ix\cdot\xi}).$$
A classical pseudo-differential operator of order $m$ has a full symbol which
is an asymptotic sum of terms $p_{m-j}(x,\xi)$ which are smooth in $\xi\ne 0$
and for $\lambda>0$ are homogeneous of degree $m-j$
$$p_{m-j}(x,\lambda\xi)=\lambda^{m-j}p_{m-j}(x,\xi).$$
The principal symbol is $p_m$ also denoted by $\sigma_m(P)$. The class of
classical pseudo-differential operators is denoted by
$\psDOcl{m}(\mathbb{R}^n)$. There are more general classes of
pseudo-differential operators based on more general symbols, but we shall
not need them here. These classes form a graded algebra under composition
$$\circ:\psDOcl{m} \times \psDOcl{m'} \rightarrow 
\psDOcl{m+m'}.$$ To obtain the principal symbol of the composite
one just takes the product: $$\sigma_{m+m'}(PQ)=
\sigma_{m}(P)\sigma_{m'}(Q)$$  however the full symbol of the product is
rather more complicated. 
Operators in $\psDOcl{-\infty}=\bigcap\limits_{m\in \mathbb{R}}\psDOcl{m}$ are
called smoothing operators. The full symbol of a pseudo-differential operator
determines the operator modulo smoothing operators. For brief introduction to
pseudo-differential operators we recommend the notes~\cite{jnotes} and for
more detail Shubin~\cite{shubin}. We note that the definition of
pseudo-differential operators can be extended to smooth manifolds using
coordinate charts.  Here the principal symbol is invariantly defined 
as a function on the cotangent bundle while the
total symbol depends on choice of coordinates.

Following~\cite{jmcd} we will consider pseudo-differential operators on a smooth manifold $Y$ depending 
smoothly on a parameter $t$.  For our purposes we will have $Y=\p M$ and $t$ the normal distance from 
the boundary.  We say that $P \in \psDO^{m,r}(Y,\mathbb{R}^+)$ if it is a family of pseudo-differential 
operators of order $m$ on $Y$, varying smoothly up to $t=0$, and such that
$$P= \sum\limits_{j=0}^{r} t^{r-j}P_j$$
with $P_j$ a smooth family of operators on $Y$ of order $m-j$. This definition extends naturally to 
operators on bundles, in our case the bundle of $k$-forms being the important example.

The symbol of $P \in \psDO^{m,r}(Y,\mathbb{R}^+)$  is defined to be the vector
$$\left( \sigma_{m-j}(P_j)\right)_{j=0}^r$$
evaluated at $t=0$.  This is a vector of functions on the cotangent bundle of $Y$.  For the case of an operator on a vector bundle, each of these functions is a field of enodmorphisms on the fibres of the bundle.

Let $u$ be a 1-form then the Bochner Laplacian (sometimes called the
rough Laplacian) is the operator expressed  in coordinates as $
-\sum\limits_{ij}g^{ij} u_{k;ij}$. The principal symbol in this case
is $g I$ where $I$ is the identity on 1-forms. 
The
formal adjoint with respect to a metric of the exterior derivative on
$k$-forms is $\delta = (-1)^{(nk+n+1)}*d*$.  The Laplacian on
differential forms (or Hodge Laplacian) is then $\Delta = d \delta +\delta d$. The principal symbol
of $d$ is $\sigma^1_d(\xi) \omega = i\xi \wedge \omega$. Let $
\contr{X}{\omega}$ denote the contraction of the differential form
$\omega$ with respect to the vector field $X$. We denote by
$\xi^\sharp$ the vector field dual to the one form $\xi$. The
principal symbol of $\delta$ is then $\sigma^1_\delta(\xi)
=-i\contr{\xi^\sharp}{\,}$. We conclude using that contraction is an
antiderivation on forms $$\sigma^2_{\Delta}(\xi)\omega =
\contr{\xi^\sharp}{(\xi \wedge \omega)} + \xi
\wedge(\contr{\xi^\sharp}){\omega} = g(\xi,\xi)\omega.$$ The
connection between the Laplacian and the Bochner Laplacian, as well as
an alternative way to calculate the principal symbol of the former, is
given by the coordinate expression for the Laplacian
\begin{eqnarray*} (\Delta u)_{i_1 \dots i_k}
&=&
\sum\limits_{ij}\left(-g^{ij} u_{i_1 \dots i_k;ij}+\sum\limits^k_{
\alpha=1}R^j_{i_\alpha} u_{i_1 \dots i_{\alpha -1}ji_{\alpha + 1}\cdots
i_k}\right.\\
& & { }\left. + \frac{1}{2} \sum\limits^k_{ \alpha=1}
\sum\limits^k_{\beta =1}R^{ij}{}_{i_{\beta} i_{\alpha}}u_{i_1 \dots
i_{\alpha -1}ji_{\alpha + 1}\dots i _{\beta -1}ii_{\beta + 1}\dots
i_k}\right). \end{eqnarray*}
 Note that for a flat
metric the Laplacian and Bochner Laplacian coincide.
A differential form $u$ satisfying Laplace's equation $\Delta u=0$ is called a
harmonic form. On a compact manifold without boundary, this is equivalent to
the condition that the form is a harmonic field, that is it is both exact,
$du=0$, and co-exact, 
$\delta u=0$ as 
\begin{equation}\label{eqn:integ}\la u,\lap u \ra = \vert\vert du
\vert\vert^2 + \vert\vert \delta u \vert\vert^2 + \int_{\p M}\delta u\wedge
*u+\int_{\p M} u\wedge*du
\end{equation}
However on manifolds with boundary there can be harmonic forms which
are not harmonic fields.

Closely related  systems of elliptic
partial differential equations occur in electro-magnetics (the vector
Helmholtz equation) and in linear elasticity. 

In a  linear elastic solid with metric tensor $g$ and with no body
forces, the displacement field $u$  (a vector field) satisfies
the equation $ {\rm Div} (CL_u g)=0$ where $u \mapsto L_u g$ is the Lie
derivative of the metric which in components is $ (L_u g)= u_{i;j}+
u_{j;i}$ (as usual a semi-colon indicates covariant differentiation
with respect to following indices) and $\rm Div$ is its formal adjoint
$a_{ij} \mapsto \sum\limits_{jk}a_{ij;k}g^{jk}$. (All sums will be
indicated explicitly.) The elastic tensor $C$ is a field of
automorphisms of the symmetric tensors on each fibre. The principal
symbol of the elastic operator is $C$. For an isotropic solid $C =
\lambda g\otimes g^\sharp + \mu I$ where $I$ is the identity operator
on symmetric tensor fields. The problem considered by \cite{nu} was
the recovery of the Lam\'{e}  parameters $\lambda$ and $\mu$ for an
isotropic solid. They also considered a related anisotropic problem
for two-dimensional elastic media.

Nakamura and Uhlmann~\cite{nu2} derive a factorisation for the anisotropic 
linear elasticity operator in boundary normal coordinates (for the flat metric).
This allows them to recover the full Taylor series of the 
`surface impedance tensor', which is a function of $C$, but not the complete 
Taylor series of $C$. For a special case, transversely isotropic media, $C$ can be 
recovered~\cite{NTU}.

In electro-magnetic theory the electric field $E$ and magnetic filed $H$
are naturally defined as 1-forms, as to take measurements of these fields one must integrate
over curves. The resulting electric and magnetic fluxes, $D$ and $B$ are naturally two forms
as one must integrate them over surfaces to make a measurement. The material properties (for simplicity 
we consider a non-chiral, linear,  insulating  material) are the permittivity $\varepsilon$ and 
permeability $\mu$, these map one forms to two forms and are the Hodge star operators for an associated electric and 
magnetic Riemannian metric. Assuming all fields to be time harmonic with angular frequency $\omega$
and the electric charge density to be constant
we have Maxwell's equations
\begin{eqnarray}\label{eqn:Max}
dB &=0,\,dD=0\\
dE &= -i\omega \mu H,\, D= \varepsilon E \\
dH &= i \omega \varepsilon E,\,B= \mu H
\end{eqnarray}
In a physically artificial situation where $\mu=\varepsilon=*$ 
(obviously after units have been scaled) we notice that $E$ and $H$ salsify the
vector Helmholtz equations
$\Delta E= \omega^2 E$ and $\Delta H = \omega^2 H$.

The main result we prove is an extension to $k$-forms for $k>0$, of
the result of Lee and Uhlmann~\cite{lu}  using a similar factorization in boundary 
normal coordinates.
Initially, our notion of a Dirichlet-to-Neumann map is non-standard for $k\ne
0$. Employing  the multi-index notation $I=(i_1,\dots,i_k)$ we write a
$k$-form 
as $u = \sum\limits_{I}u_I dx_I$ where $dx_I=dx_{i_1}\wedge
dx_{i_2}\wedge\dots \wedge dx_{i_k}$. Following~\cite{lu} we use a coordinate chart
$(x_1,\dots,x_n)=(x',x_n)$ where $x'=(x_1,\dots,x_{n-1})$ is a chart on the
boundary, and $x_n$ is the distance to the boundary. We denote by
$\p_n$ both a unit vector field normal to the boundary and its
associated derivation on functions. We extend this to $k$ forms by its
actions on the components as functions
$\p_n u =  \sum\limits_{I}\p_n u_I\, dx_I$. The operator
$\Lambda_g:u|_{\p M} \mapsto (\p_n u)|_{\p M}$, where $\Delta  u =0$, is
linear while somewhat unnatural. We will discuss the relationship to
more natural Neumann data for Laplace's equation in Section 2. Our
main result is 

\begin{thm}\label{prop:mainthm}
Let $M$ be a smooth compact orientable Riemannian manifold with
boundary, with $\mathrm{dim}(M)=n>2$ and metric $g$, and let $k$ be an
integer $0\le k \le n$. 

\noindent (i) The Dirichlet-to-Neumann map  $\Lambda_g:u|_{\p M}
 \mapsto (\p_n u)|_{\p M}$ for the $k$-form Laplace's
equation $\Delta u =0$ is a classical pseudo-differential operator of
order one.

\noindent(ii) The Taylor series, at the boundary, of the metric  in boundary normal coordinates is uniquely determined by
the full symbol of $\Lambda_g$.  For $0<k<n$  only one 
diagonal component of the full symbol is needed corresponding to
$dx_I=dx_{i_1}\wedge dx_{i_2}\wedge\dots \wedge dx_{i_k}$
but for  $k=(n+1)/2$ the multi-index $I=(i_1,\dots,i_k)$ must
exclude $n$ and for $k=(n-1)/2$, $I$ must include $n$. \end{thm}

Where our work 
differs from Lee and Uhlmann's is in its use of families of 
operators parameterized by the normal
distance. The case $k=n$ is clearly equivalent to $k=0$ so we need only
consider the case $0<k<n$. 
Lee and Uhlmann showed that the full Taylor series of a real analytic
metric on a real analytic manifold, where the relative homotopy group
of the boundary $\pi(M,\partial M)$ is trivial, determines the metric,
provided the manifold is strongly convex or the metric can be extended
analytically to a larger manifold without boundary. Recent work of
Lassas and Uhlmann~\cite{lau} removes these geometric and topological
hypotheses showing that an analytic metric is determined throughout
any connected analytic manifold by the Dirichlet-to-Neumann map. It remains to be seen if the same techniques can be
applied to the general $k$-form case.

Only a small modification of the argument is needed to prove a version
of Theorem~\ref{prop:mainthm} for the equivalent $k$-form Helmholtz problem at fixed
frequency. %
\section{Boundary Conditions}\label{section:boundcond}
%
%
In a neighbourhood of the
boundary where our boundary normal coordinates are defined, we can
distinguish  between {\em tangential} $k$-forms which have no $dx_n$
in their coordinate expression, and {\em normal} forms which have a
common factor $dx_n$. Clearly the space of $k$-forms on this
neighbourhood  is the direct sum of the spaces of tangential and
normal forms. The projection on to the tangential component is
$\pi_{\mathrm t}(\omega) = \contr{\p_n}{(dx_n\wedge \omega)}$ and on to the
normal component  $\pi_{\mathrm n}(\omega) = dx_n\wedge
(\contr{\p_n}{\omega})$. From these expressions it is clear that $*
\pi_{\mathrm n} = \pi_{\mathrm t} *$ and that the Hodge star of a tangential form is normal
and {\em vice-versa}.

Let $i:\partial M
\rightarrow M$ be the inclusion of the boundary. 
The tangential
component of $k$-form $u \in \Omega^k(M)$ at the boundary is then 
 the pull-pack to the
boundary $i^*u  \in \Omega^k(\partial M)$.
  The normal part of $u$ at the boundary, can be
determined uniquely from $i^**u = *_\p(\contr{\p_n}{u})|_{\p M} \in
\Omega^{n-k}(\p M)$ where $*_{\p}$ is the induced Hodge-star on the
boundary. Dirichlet data for harmonic
$k$-forms    
consists of the both the tangential  component  and the normal
component 
of the form at 
the boundary~\cite{ds}. Note a possible source of confusion here. When considering
forms on manifolds with a boundary $i^*u$ is often thought of as the restriction of form to the
boundary, and the pull back is often omitted, for example in Stoke's
theorem. In Theorem~\ref{prop:mainthm} our Dirichlet data $u|_{\p M}$
is the restriction to the boundary of a form, but in the sense of considering only base
points on $\p M$ and no projection of the fibre on to the tangential component. This Dirichlet data
together with the integrals of $u$ on a basis of the relative homology group $H_k(M,\p M)$
determines a unique solution to $\Delta u =0$. For simplicity we will assume that the said integrals are 
specified to be zero. We note that for the case $k=0$ one simply specifies the integral of $u$ on each 
connected component of $M$ with a non-empty boundary. The natural Neumann data~\cite{ds} is the 
specification of $i^**du$ and $i^*\delta u$. Note that as in the case $k=0$ where Neumann data $i^**du$ must have
zero integral on the boundary, there are compatibility conditions for
Neumann data~\cite{ds}. A  natural Dirichlet-to-Neumann mapping is
therefore $\Pi_g: \Omega^k(\p M)\times\Omega^{n-k}(\p M) \rightarrow
\Omega^{n-k-1}(\p M)\times\Omega^{k-1}(\p M)$ given by
$$ \Pi_g(f_\tau,f_\nu)= (i^**du,i^*\delta u)$$
where $\Delta u = 0$, $i^*u=f_\tau$, $i^**u=f_\nu$. Here we use $\nu$
and $\tau$ as labels for the normal and tangential components in the
sense of the domain and range of $\Pi_g$. We can now recast Theorem
\ref{prop:mainthm} in terms of this natural data.

\begin{cor}\label{cor:natural}
Let $M$ be a smooth compact orientable Riemannian manifold with
boundary, with $\mathrm{dim}(M)=m>2$ and metric $g$, and let $k$ be an
integer $0\le k \le n$.  
\noindent 
(i) The natural
Dirichlet-to-Neumann map  $\Pi_g$ for the $k$-form Laplace's equation defined
above  is a classical pseudo-differential operator of
order 1.

\noindent 
(ii) The full symbol of $\Pi_g$ determines
the Taylor series (at the boundary) of the metric in boundary normal
coordinates.  Furthermore for    $k\not\in\left\{(n-2)/2, (n-1)/2,n\right\}$ only the full symbol of the
tangential part $\Pi_{g\, \tau\tau}$ is needed and for  and  $k\not\in
\left\{0, (n+1)/2,(n+2)/2\right\}$ only
the full symbol of the normal part $\Pi_{g\, \nu\nu}$ is needed.
\end{cor}
\noindent The proof of Corollary \ref{cor:natural} follows the proof of Theorem
\ref{prop:mainthm} in Section~\ref{sec:factsymb}.
 
In the case of electromagnetics note that when formulated in differential forms 
Maxwell's equations~(\ref{eqn:Max})--(1.4) are independent of the ambient Euclidean metric and
thus invariantly defined. In the inverse boundary value problems for time 
harmonic Maxwell's equations one typically specifies the boundary data
invariantly as the tangential component of $E$ and $H$.
The isotropic case where both the electric and magnetic metrics are conformally flat,
has been studied by~\cite{OPS} and~\cite{jmcd}. 

By contrast in elasticity strain is a measure of the distortion of the Euclidean metric,
and one seeks the elastic tensor with the ambient metric given. This problem is not
diffeomorphism invariant.
\section{Factorization and symbol calculation}\label{sec:factsymb}

We consider metrics to be equivalent if they are related by
a diffeomorphism which fixes  points on the boundary.
Without loss of generality, therefore,  we  can assume 
that $x_n$ is the boundary
normal coordinate for both metrics. Later we will make a more
specific choice for the coordinate chart on the boundary.

We use notation from \cite{jmcd}, in particular
$\psDO^{m,l}$ denotes families of pseudo-differential operators,
$P_{x_n},$ in $x'$ such that the $j$ term in the total symbol vanishes
to order $l-j$ at $x_n=0,$ and $\mathrm{DO}^{m,r}$ is the class of
such differential operators.

Let $\Delta' = \sum \limits_{i,j=1}^{n-1}
h^{ij}D_{x_i}D_{x_j}.$ We have that $$\Delta = (D_{x_n}^{2} + \Delta')
I + ED_{x_n} + H(x,D_{x'}),$$ where $H$ is a first order system in
$x'$  on the bundle of $k$-forms and $E$ is
a smooth endomorphism of the bundle of $k$-forms. We use the notation
$|\xi|= \sqrt{g(\xi,\xi)}$ for a covector $\xi$.
\begin{prop} \label{prop:factorization}
There exists a $B(x,D_{x'}) \in \psDOcl{1}(\p M; \mr ;
\Omega^{1}(M))$ such that
$\sigma_{1}(B) = |\xi'|_{x}I$ and
$$ \Delta = (D_{x_n}I + E + iB)(D_{x_n}I - iB),$$
modulo smoothing
and $B$ is unique modulo smoothing.
\end{prop}
\begin{proof}
If we expand, we obtain
$$ D_{x_n}^{2} I + B^2 + ED_{x_n}-iEB+ i[D_{x_n},B].$$
Taking the principal symbol of $B$ as specified we obtain an error, in
$\psDOcl{1}.$  Now suppose we have chosen $B_{j}$ such that the error,
$F_{j},$ is in $\psDOcl{1-j}.$ Let $B_{j+1} = B_{j} + C$ with $C$ in
$\psDOcl{-j}.$ Upon expanding we then obtain an extra term $CB_{j} + B_{j}C +
E_j$ with $E_{j}$ of order $-j.$ Taking $\sigma_{-j}(C) =
-\frac{1}{2}|\xi'|_{x}^{-1}\sigma_{1-j}(F_j).$ We have achieved an error one
order better. Inducting and summing, we achieve an error in $
\psDOcl{-\infty}.$

As the choice at each stage was forced, $B$ is unique.
\end{proof}

The importance of this factorization is that $B(0,D_{x'})$ is equal
modulo smoothing terms to $\Lambda_g$. We will summarise the
argument which is identical to that given by~\cite{lu} for the
$0$-form case. Given a harmonic $k$-form $u$, let $v= (D_{x_n}I -
iB)u$ so that $(D_{x_n}I + E + iB)v=0$. These are both generalised
heat equations, the second with `time' reversed. As both are
smoothing we see that $\partial_n u  = Bu + \mbox{smooth terms}$ and hence
$\Lambda_g=B \mod \psDOcl{-\infty}$.  This proves
Theorem~\ref{prop:mainthm} part (i) and part (ii) will
follow if we can show that two metrics $g_1,g_2$ with identical full
symbols of $B$ at the boundary must agree to infinite order at the
boundary. Rather than calculating the full symbol of $B$, we use the
calculus of pseudo-differential operators parameterised by the normal
distance. The advantage is that we need only calculate principal
symbols.

We want to compare the Laplacians
associated to two different metrics assumed equal up to order $l$ in
the normal coordinate. Of course it is immediate from the principal
symbol of $B$ that the metrics agree on the boundary so we can take
$g_{1} - g_{2} = x_{n}^{l} k$ for some $l>0$

First we compare the Hodge star operators.
By definition $ \omega \wedge *\omega = g( \omega,\omega)
\mu$ and we see that
$$*_1 = *_2 + x_{n}^{l}\alpha$$
where $\alpha$ is a smooth homomorphism from
$\Omega^{k}$ to $\Omega^{n-k}.$

\begin{lemma}
If $\Delta_{j}$ is the Laplacian on  $k$-forms associated with $g_j$ then
$$\Delta_{2} - \Delta_{1} = x_{n}^{l-1}FD_{x_n} + A$$ where $F$ is a smooth
endomorphism and $A \in \mathrm{DO}^{2,l}.$
\end{lemma}
\begin{proof}
The Laplacian is defined by $*d*d$ and $d*d*$ where $d$ is independent of the
metric and $*_2 = *_1 + x_{n}^{l}\alpha.$

The result follows simply by observing that in
$d *_2 d *_2$ and $*_2 d *_2 d$ terms not in $\Delta_1$ will vanish to order
$l$ at $x_n=0$ unless $d$ is applied to the $x_{n}^{l}$ term. If $d$ is applied
once to such a term we obtain a first order differential operator vanishing to
order $l-1$ and if twice a zeroth order operator vanishing to order
$l-2.$ This is the result desired --- we know there are no second
order terms in $D_{x_n}$ from our expression for the principal symbol.
\end{proof}

\begin{lemma}
Let $\Delta_{j}$ be factored as in Prop \ref{prop:factorization} with $E_j,
B_j$ the corresponding terms. We then have that,
$$B_2 - B_1 \in \psDO^{1,l}.$$
\end{lemma}
\begin{proof}
Let $C=B_2 - B_1 .$ As the
principal symbols of $B_1, B_2$ agree at $x_n=0$ so we have that $C$ is in
$\psDO^{1,1}.$

Note that $E_2 = E_1 +x_{n}^{l-1}F.$ Expanding the factorizations for
$\Delta_2, \Delta_1$ and subtracting we have that,
$$\Delta_2 - \Delta_1 = i[C,D_{x_n}] + B_1 C +CB_1 + x_{n}^{l-1}FD_{x_n} -
x_{n}^{l-1}iF(B+C)
+C^{2}.$$ After cancelling the $x_{n}^{l-1}FD_{x_n},$ we have that
$$i[C,D_{x_n}] + B_1 C +CB_1 -ix_{n}^{l-1}F(B+C) +C^{2} \in \psDO^{2,l}.$$
If $C \in \psDO^{1,r}$ with $1 \leq r < l$ then we have that $[C,D_{x_n}] \in
\psDO^{1,r-1},$ $B_1 C +CB_1 \in \psDO^{2,r},$ $C^{2} \in \psDO^{2,2r},$
 $x_{n}^{l-1}iFB \in \psDO^{1,l-1},$ and $x_{n}^{l-1}iFC \in 
\psDO^{1,r+l-1}.$
Taking the residue modulo $\psDO^{2,r}$ we have that
$i[C,D_{x_n}] + B_1 C+ CB_1$ is congruent to zero modulo $\psDO^{2,r}.$
Recall that
$$C = \sum \limits_{j \leq r} x_{n}^{r-j}C_j$$
with $C_j \in \psDOcl{1-j}.$ The only term of second order is $C_1 B
+BC_1$ so we deduce that the principal symbol of $C_1$ vanishes at
$x_n=0.$ Let $c_j$ denote the principal symbol of $C_j$ at $x_n=0.$ We
then have that $$(r-j) c_j + 2|\xi'|_{x}c_{j-1}=0$$ for each $j$.
Iterating we conclude that $c_j=0$ for each $j$ which
proves that $C \in \psDO^{1,r+1}.$ Repeating, the result follows.
 \end{proof}

It follows from this Lemma that $B_2-B_1$ restricted to $x_n =0$ is $C_l$ a
pseudo-differential operator of order $1-l.$ Our main result will follow if we
can compute the principal symbol of this operator and show that it
determines $k$, the lead term of $g_1 - g_2.$


Now
$$\Delta_2 - \Delta_1 = x_{n}^{l} P_2 + x_{n}^{l-1} P_1
+x_{n}^{l-2}P_{0} +x_{n}^{l-1}FD_{x_n} $$
with $P_j$ a differential operator in $x'$ of order $j$
and we know from our principal symbol computation that
$P_2$ is equal to
$\sum\tilde{k}_{ij}D_{x_i}D_{x_j}$ where $\tilde{k} = - h^{-1}kh^{-1}$ where
$g =dx^{2} + h(x',dx') + \mathcal{O}(x_n).$

Arguing as above with $C=B_2 - B_1$ we have
\begin{equation}
i[C,D_{x_n}] + B_1 C +CB_1 -ix_{n}^{l-1}F(B+C) +C^{2} = x_{n}^{l} P_2
+ x_{n}^{l-1} P_1 +x_{n}^{l-2}P_{0}.
\end{equation}
So modulo $\psDO^{2,l}$ we have,
\begin{equation}
 i[C,D_{x_n}] + B_1 C +CB_1 -ix_{n}^{l-1}FB_1 = x_{n}^{l} P_2 +
x_{n}^{l-1} P_1 +x_{n}^{l-2}P_{0}.
\end{equation}
Let $C= \sum \limits_{j=0}^{l} x_{n}^{l-j}C_{j}$ with $C_j \in
\psDOcl{1-j},$ let $c_j$ be the principal symbol of $C_j$. We have,
\begin{align*}
2|\xi'|_{x} c_{0} =& \sum\limits_{i,j<n} \tilde{k}_{ij} \xi_i \xi_j,\\
2|\xi'|_{x} c_{1} + lc_{0} =& \sigma_{1}(P_1)(\xi') + i|\xi'|_{x} F,
\\ 2|\xi'|_{x} c_{2} + (l-1)c_{1} =& \sigma_{0}(P_0), \\
2|\xi'|_{x} c_{2+j} + (l-j-1)c_{1+j}=& 0, \text{ for } 1<j \leq l-2.
\end{align*}
We therefore deduce that
\begin{equation}\label{equation:crelations}
c_l = K_{l} (|\xi'|_{x})^{-l-1} \sum\limits_{i,j<n} \tilde{k}_{ij}
\xi_i \xi_j I +L_{l} |\xi'|^{-l} \sigma_{1}(P_1)(\xi') + iL_{l}
|\xi'|_{x}^{1-l} F + M_l \sigma_{0}(P_0)|\xi'|_x^{1-l}
\end{equation}
where $K_l, L_l, M_l$ are computable non-zero constants. (Of course,
$M_0, M_1 =0.$ )

We want to show that $c_l=0$ implies that $\tilde{k}_{ij}=0$ for all $i,j.$

If
the principal symbol $c_l=0,$ we have taking any component $rr$ of the symbol that
$$K_{l} \sum \limits_{ij} \sum \tilde{k}_{ij}(x) \xi_i \xi_j
+L_{l} |\xi'|\sigma_{1}(P_1)_{rr}(x,\xi') + |\xi'|^{2}(iL_{l}F_{rr} +
C_l \sigma_{0}(P_{0})_{rr})=0.$$
As $P_1$ is a differential operator $\sigma_{1}(P_1)_{rr}(x,\xi')$ is
linear in $\xi'$ and the final two terms are independent of $\xi'.$
Since the middle term is not smooth, unless zero, as $\xi' \to 0$ and
the other terms are smooth, we deduce that $\sigma_1(P_1)=0$ and
so we have,
$$K_{l} \sum \limits_{i,j<n} \sum \tilde{k}_{ij}(x) \xi_i \xi_j+
|\xi'|^{2}(iL_{l}F_{rr} + C_l \sigma_{0}(P_{0})_{rr})=0.$$
This shows that $k_{ij}$ must be a scalar multiple of the
identity matrix. Or more invariantly that $k(x)$ is a scalar multiple
of $h(x).$ To see that $k$ is actually zero, we need to compute more
precisely.

For convenience we now
reduce to the case of a Euclidean background metric. We first prove
\begin{lemma}Let \begin{align*}g_1 &= h + x_n m, \\g_2 &= h+ x_n m +
x_{n}^{l} r,\end{align*} and let $*_j$ be the Hodge $*$ operator
associated to $g_j.$ We then have that $*_1 - *_2$ modulo ${\mathcal
O}(x_{n}^{l+1})$ is independent of $m.$ \end{lemma}
\begin{proof}Let $\mu$ be a volume form
of $h$ and $\mu_j$ a volume form for $g_j$. We have by definition
$h(\mu,\mu)=1$ and  $g_j(\mu_j,\mu_j)=1$.
 We then have that $\mu_2=(1+x_n m(\mu,\mu) +x_{n}^{l}
r(\mu,\mu))^{-1/2}\mu$ and similarly
for $g_1.$ It is now clear that

$$\mu_2= \frac{1}{\sqrt{ 1+x_n m(\mu,\mu)+{x_l}^l r(\mu,\mu)
}} \mu. $$
and
$$\mu_1= \frac{1}{\sqrt{ 1+x_n m(\mu,\mu)}} \mu $$
For any $\nu,\omega \in \Omega^k(M)$

\begin{eqnarray*}
\nu \wedge (*_2 - *_1) \omega &=& \frac{(h+x_n m +
x_{n}^{l}r)(\nu,\omega)}{\sqrt{1+x_n m(\mu,\mu) + x_{n}^{l}
r(\mu,\mu)}} - \frac{(h +x_{n} m)(\nu,\omega)}{\sqrt{1+
x_{n}m(\mu,\mu)}} \\
&=& \left( {x_n}^l \left( r(\nu,\omega) -
\frac{1}{2}r(\mu,\mu)h \right) + {\mathcal O}({x_n}^{l+1})\right)\mu
\end{eqnarray*}
which does not involve $m.$\end{proof}
We also have,
\begin{lemma}\label{prop:euclmet}Let $g_2 = dx^2
+ x_n^{l} r + m,$ and $g_1 = dx^2 + m,$ where $m$ vanishes to second
order at the origin. Let $*_j$ be the Hodge $*$ operator of $g_j.$
We then have that $*_2 - *_1$ is independent of $m$ modulo terms of
the form $x_{n}^{l} t + x_{n}^{l+1} w$ where $t$ vanishes to second
order at the origin and $w$
is smooth. \end{lemma}
\begin{proof}Let $\mu$ be
the volume form for $dx^2.$ Arguing as above we have
that,
\begin{equation}\nu \wedge (*_2-*_1) \omega = \frac{(dx^2 + m +
x_{n}^{l} r)(\nu,\omega)}{\sqrt{1+m(\mu,\mu) +x_{n}^{l}r(\mu,\mu)}} -
\frac{(dx^2 +m)(\nu,\omega)}{\sqrt{1+m(\mu,\mu)}},
\end{equation}
which modulo
terms vanishing appropriately at $x=0$ equals the bilinear
form
\begin{equation}(dx^2 + m + x_{n}^{l}r)(1-
\frac{1}{2}(m(\mu,\mu)+x_{n}^{l}r(\mu,\mu)) - (dx^2
+m)(1- \frac{1}{2}m(\mu,\mu)),
\end{equation}
which upon expanding modulo appropriately vanishing
terms  equals,
\begin{equation}
x_{n}^{l} r - \frac{1}{2}x_{n}^{l} r(\mu,\mu)
dx^2,\end{equation}
which does not involve $m.$ \end{proof}
Now fix a point
$p$ where we will calculate the principal symbol of the difference of
the Dirichlet to Neumann maps and show that it being zero implies that
the next term of the difference of the metrics also vanishes there. We
take geodesic normal coordinates about $p$ in the boundary and then
extend normally with respect to $g_1.$ 
Now we fix  a point on the boundary $p$ where we will calculate the 
principal symbol of the difference of
the Dirichlet to Neumann maps and show that it being zero implies that
the next term of the difference of the metrics also vanishes there.
We choose $x'$ to be a
Riemann normal coordinate system on the boundary and then
extend normally with respect to $g_1.$  In particular, we
have $$g =dx_{n}^{2} + h(x,dx'),$$
and on the boundary $h((x',0),dx') = d{x'}^2 +\mathcal{O}({x'}^2)$.

We then have, using
Lemma~\ref{prop:euclmet} that,
$$ g_1 = dx^{2} + m ,$$ and $$g_2 = dx^2 + m +
x_{n}^{l}r$$ where $m= x_{n}^{l} t$ with $t$ vanishing to second order at
$x'=0.$ Let $g_{1}^{'} = dx^{2},$ and $g_{2}^{'} = dx^{2} + x_{n}^{l} r.$ We
then have that $*_{j} - *_{j}^{'}= x_{n} \alpha_{j} + \beta_{j}$ for $j=1,2$
with $\beta_j$ vanishing at $x'=0$ to second order. We also have from our
lemmas that $$ (*_1- *_2) - (*_{1}^{'} - *_{2}^{'}) = x_{n}^{l}
\gamma$$ with $\gamma$ vanishing to second order at $x'=0.$

It is now
clear that when computing the lead term of $d*_2 d (*_2 -*_1)$ and its
appropriate permutations at the point $p$ that we can replace $*_j$
by $*_{j}^{'},$ without  changing the value. So to finish our theorem
we take $g_1 = dx^{2}$ and $g_2 = dx_{n}^2+(1+x_{n}^{l}\lambda(x'))
d{x'}^{2},$ with $\lambda$ a smooth function.

Let us consider the action of $*_2$ on normal and tangential forms.
For a multi index $I=(i_1,i_2,\dots,i_k)$ we use the convention
$u_Idx_I=u_{i_1i_2\dots i_k}dx_{i_1}\wedge dx_{i_2}\wedge\dots
dx_{i_k}$. We will denote by $I'$ the complimentary multi index with
$I'=(j_1,\dots,j_{n-k})$ where
$(i_1,\dots,i_k,j_1,\dots,j_{n-k})$ is an even permutation of
$1,\dots,n$.

For the metric $g_1$,$(\p_1,\dots,\p_n)$ is an oriented
orthonormal frame and $(dx_1,\dots,dx_n)$ an oriented orthonormal
co-frame. We get an orthonormal frame for $g_2$ by dividing each
$\p_i$ by $(1+\lambda(x')x_{n}^{l})^{1/2}$ except for $\p_n$.
Similarly $\eta_i = (1+\lambda(x')x_{n}^{l})^{1/2}dx_i$ for $i<n$
and $\eta_n=dx_n$ is an orthonormal coframe for $g_2$. Applying $*_2$
to a typical normal basis $k$-form  $dx_I$ (with $n\in I$), we have
\begin{align*}
*_2(dx_I)=&
(1+\lambda(x')x_{n}^{l})^{-\frac{k-1}{2}}*_2(\eta_I)\\
=& (1+\lambda(x')x_{n}^{l})^{-\frac{k-1}{2}}(\eta_{I'})\\
=& (1+\lambda(x')x_{n}^{l})^{\frac{n+1}{2}-k}(dx_{I'}).
\end{align*}
For a typical basis element of $\Omega^k_t(\p M)$, we have $n\not\in
I$ and
\begin{align*}
  *_2(dx_I)=&
(1+\lambda(x')x_{n}^{l})^{-\frac{k}{2}}*_2(\eta_I)\\
=& (1+\lambda(x')x_{n}^{l})^{-\frac{k}{2}}(\eta_{I'})\\
=& (1+\lambda(x')x_{n}^{l})^{\frac{n-1}{2}-k}(dx_{I'}).
\end{align*}
Modulo $x_{n}^{l+1},$ we have
that on normal $k$-forms $(*_2 - *_1) = x_{n}^{l} \left(
\frac{n+1}{2}-k\right)\lambda *_1$  and on tangential  
$k$-forms $(*_2 - *_1) = x_{n}^{l} \left(\frac{n-1}{2}-k\right)\lambda *_1.$ 

To prove the theorem we only have to consider the action of
the difference of the operators on a particular $k$-form $u_Idx_I$.

 We first establish that there is no
contribution from $$*_2 d *_2 d - *_1 d*_1 d= (*_2-*_1)d *_2 d +*_1
d(*_2 - *_1)d.$$ Any first order term in the first term of the RHS
will be in $\mathrm{DO}^{1,l}$ and therefore not contribute to the
principal symbol.
For the second term in the RHS, we have,
$$ *_1 d(*_2 - *_1)d (u_{I}dx_I)
=  *_1 d\left( \left({\frac{n\mp 1}{2}-k}\right)x_{n}^{l} \lambda(x')
\right. \left. *_1 \sum_{j\not\in I}\frac{\p u_I}{\p x_j} dx_j \wedge dx_I \right)
$$
plus terms involving $x_{n}^{l+1}$ which
will not contribute. Now consider the second $d,$ if it applies to
$u_I$ we get a second order term which we already
understand; if it applies to $x_n$ we get a term involving $dx_n$
which will then have no $dx_n$ component on applying $*_1$ again. Thus
we get no sufficiently low order contribution to either the zeroth
order term or to the coefficient of $D_{x_n}.$

Moving on to
$$d*_2 d *_2  - d*_1 d*_1 = d(*_2-*_1)d *_2  +d*_1d(*_2 - *_1),$$
looking at the first term on the RHS we can equally compute with
$d(*_2-*_1)d*_1$ as the difference will be in a non-contributory
residue class.

First considering the normal case $n\in I$ we compute
$$d*_1u = \sum
\limits_{j\not\in I'}\p_j u_I\, dx_{j} \wedge dx_{I'}.$$
Now consider
\begin{align*}
&d(*_2-*_1)(\p_j u_I dx_{j} \wedge dx_{I'}) \\
&=\left(\frac{n\pm 1}{2}-k\right)d\left(x_{n}^l\lambda \p_j u_I
 *_1(dx_{j} \wedge dx_{I'}) \right)\\
&= \left(\frac{n\pm 1}{2}-k\right)\left( lx_{n}^{l-1} \lambda  \p_j
u_{I} dx_n + {x_n}^{l} d(\lambda  \p_j u_I)
\right)\wedge*_1(dx_{j} \wedge dx_{I'})
\end{align*}
Where the $+$ holds for $j=n$ and the $-$ otherwise,
but for  this to have a $dx_I$ component  we
must have $j=n$. The second term in the final bracket has a
coefficient $x_n^l$ so can only contribute to the second order part
which we already understand. We do  have a contribution to
$F$ from
\begin{equation}\label{equation:bitofF}
 \left(\frac{n+1}{2}-k\right)
lx_{n}^{l-1} \lambda \p_n u_I dx_I.
\end{equation}

We are left with the contribution of $d*_1 d(*_2 - *_1)$.
If we apply $d$ to $(*_2-*_1)(u_Idx_I)$, we can drop the terms where
$d$ falls on a $\lambda$ as these are in $\psDO^{1,l}$. So on applying
$d(*_2 - *_1)$ to $ u_I dx_I$ we are left with,
$$\left( \frac{n+1}{2}-k\right)\left( x_n^l\lambda
\sum\limits_{j\in I}\p_j u_I\,dx_j\wedge dx_{I'} +(l-1)x_n^{l-1} \lambda
u_I dx_n\wedge dx_{I'} \right).$$
On applying  $*_1,$
we get
$$ \left( \frac{n+1}{2}-k\right) \left(
 x_n^l\lambda\sum\limits_{j\in I}\p_j u_I\, dx_{I_j}
+(l-1)x_n^{l-1} \lambda
u_I  dx_{I_n}
 \right)$$
where $I_j$ is simply $I$ with $j$ deleted.

We now apply $d$
again to get the final contributions. We get another contribution
to $F$ identical to Equation~\ref{equation:bitofF} so that
$$F= 2 il \left(\frac{n+1}{2}-k\right)\lambda $$
and an $x_n^{l-2}u_I$ which gives
$$\sigma_0(P_0)=  l(l-1)\left( \frac{n+1}{2}-k\right)\lambda$$
Substituting back in
to Equation~\ref{equation:crelations} we see that if
$c_l=0$ then $\tilde{k}_{ij}=0$ and  Theorem~\ref{prop:mainthm} is
proved for the case for $k\ne (n+1)/2$.  Now we consider the case of tangential data
that is $u_Idx_I$ with $n \not\in I$.  

A similar argument applies to a tangential form, with $k\ne(n-1)/2$ and
Theorem~\ref{prop:mainthm} is proved.

\noindent{\em Proof of Corollary~\ref{cor:natural}}.  First observe that given
the induced metric on the  boundary the data $(i^*u,i^**u)$ determines
$u|_{\p M}$.  In a neighbourhood's of the
boundary a  $k$-form $u$ can be expressed as
$$u = \sum\limits_{|I|=k,n\not\in I} u_Idx_I + \sum\limits_{|J|=k-1,n\not\in
J}u_{(n,J)} dx_n\wedge dx_J$$
so that  
$$\contr{\p_n}{du} = 
\sum\limits_{|I|=k,n\not\in I} \p_n  u_I\, dx_I$$
and
$$ i^**du = *_\p (\contr{\p_n}{du})|_{\p M} = *_\p \pi_{\mathrm t} \Lambda_g \left(u|_{\p M}\right)$$ 
This shows that  $\pi_\tau \Pi_g$ is a pseudo
differential operator of order 1. Notice now  that for a harmonic
$k$-form $u$, $v=*u$ is a harmonic $n-k$ form for which both the
tangential and normal parts of Dirichlet and Neumann data are
exchanged $i^*v =  i^**u$, $i^* *v = \pm i^*u$, $i^**dv =  \pm i^*\delta u$ and
$i^* \delta v = \pm i^*du $. It follows that $\pi_\nu \Pi_g$ is also a
pseudo differential operator of order 1. Thus we have proved part (i)
of Corollary~\ref{cor:natural}.  For part (ii) notice that for $k=0$ the result
is proved by~\cite{lu}, and here $\delta u =0$ identically so the
normal part of the Neumann data gives us no information. Similarly for
$k=n$ the tangential part of the Neumann data vanishes. For the case
$0<k<n$  and $k\ne (n-1)/2$ where any tangential-tangential
diagonal component of $\Lambda_g$ determines the Taylor series we
require, it is clear that the Taylor series  is determined as
long as we have  $*_\p$. We have the  principal symbol $\sigma_1(\Pi_{g\,
\tau\tau})(\xi)= *_\p |\xi|_g$  and to finish the proof for this case
we show that this determines  $*_\p$ at each point on the boundary.

Fix a point on $\p M$ and choose any multi-indices $I_0$ and $J_0$
such that $g_0(\xi,\xi):=(\sigma_1(\Pi_{g \tau\tau})(\xi)_{I_0J_0})^2$ is a non-zero quadratic function of $\xi$. Now
$g(\xi,\xi) = \alpha g_0(\xi,\xi)$  where $\alpha = 1/(*_{I_0J_0})^2$ is
to be determined. Let $*_{0 \p}$ be  the Hodge star on $k$-forms on the
boundary determined by $g_0$ then  $*_\p = \alpha^{k - (n-1)/2}*_{0\p
}$ and $(\sigma_1(\Pi_{g\tau\tau})(\xi)_{I_0J_0} = \alpha^{k-
  (n-2)/2}*_{0 \p} |\xi|_{g_0}$. As $g_0$  is known $\alpha$ is determined
provided $k\ne (n-2)/2$, hence we have  $g$ at the boundary and $*_\p$.

For the case $k \not\in\left\{0,  (n+1)/2,(n+2)/2 \right\}$ we simply  apply the above argument to
boundary data for $*u$. As one of the  conditions on $k$ must hold it is
certainly true that the full symbol for the complete $\Pi_g$
determines the Taylor series, at the boundary,  of $g$.


\begin{thebibliography}{99}
\bibitem{ds} G.F.D. Duff and D.C. Spencer, Harmonic tensors on
Riemannian manifolds with boundary, Ann. of Math, 56, 128--156, 1952
\bibitem{gilk}
P.B. Gilkey, Invariance theory, the heat equation and the Atiya-Singer
index theorem, Publish or Perish, 1984, electronic reprint
{\tt{http://www.mi.sanu.ac.yu/EMIS/monographs/gilkey/index.html}}.
\bibitem{jnotes}M.S. Joshi, Introduction to
Pseudo-differential operators, arXiv.org e-print math.AP/9906155, 1999, 
{\tt{http://arXiv.org/abs/math.AP/9906155}}.  
\bibitem{jmcd} M.S. Joshi, S. McDowall, Total determination of
material parameters from electromagnetic boundary information, Pacific
J. Math, 193, 107--129, 2000
\bibitem{lau} M. Lassas, G.  Uhlmann, On determining a Riemannian
manifold from the Dirichlet-to-Neumann map,
 Ann. Sci. \'{E}cole Norm. Sup.,  34, 771--787, 2001
\bibitem{lu} J. Lee, G. Uhlmann, Determining anisotropic
real-analytic conductivities by boundary measurements, Comm. Pure
Appl. Math. 42, 1097--1112, 1989
\bibitem{nu}G. Nakamura and G. Uhlmann, Inverse problems at the
boundary for an elastic medium, SIAM J. Math. Anal,
26, 263--279, 1995
\bibitem{nu2}G. Nakamura and G. Uhlmann, A Layer Stripping Algorithm
in Elastic Impedance Tomography, 
in {\em Inverse Problems in Wave Propagation}, IMA Vol. Math. Appl. 90,
G. Chavent {\em et al} eds.,
Springer-Verlag, New York,  375--384, 1997. 
\bibitem{NTU} G. Nakamura, K. Tanuma and G. Uhlmann, Layer stripping for a transversely 
isotropic elastic medium, SIAM J. Appl. Math., 59,  1879--1891, 1999
\bibitem{shubin} M.A. Shubin, Pseudo-differential Operators and Spectral
Theory, Springer-Verlag, Berlin, 1985.
\bibitem{OPS} 
P. Ola, L. P\"aiv\"arinta, E. Somersalo
An inverse boundary value problem in electrodynamics. 
 Duke Math. J., 70, 617--653, 1993. 
\end{thebibliography}
 \end{document}